\documentclass[11pt]{article}
\usepackage{amsmath,amssymb,amsfonts,amsthm,dsfont,enumerate,a4wide}
\usepackage{hyperref}
\usepackage{color}

\theoremstyle{plain}
\newtheorem{theorem}{Theorem}

\theoremstyle{remark}
\newtheorem{remark}{Remark}
\newcommand{\cD}{{\cal D}}

\newcommand{\cL}{{\cal L}}

\newcommand{\cN}{{\cal N}}

\newcommand{\cS}{{\cal S}}

\newcommand{\GR}{\mathds{R}}

\newcommand{\GN}{\mathds{N}}
\newcommand{\GD}{\mathds{D}}

\newcommand{\bP}{\mathbb P}
\newcommand{\bE}{\mathbb E}

\newcommand{\conver}{\mathop{\longrightarrow}}
\newcommand{\wconver}{\mathop{\Longrightarrow}}

\newcommand{\wconv}{\wconver}
\newcommand{\conv}{\conver}

\newcommand{\indist}{\ \conver_{\cD}\ }

\newcommand{\Cov}[2]{\text{ Cov\,}(#1,#2)}

\newcommand{\bequ}[1]{\begin{equation}\label{#1}}
\newcommand{\eequ}{\end{equation}}
\newcommand{\nref}[1]{(\ref{#1})}

\begin{document}
	
	\title{Stable limits for associated regularly varying sequences}
	
	\author{Adam Jakubowski\footnote{E-mail: adjakubo@mat.umk.pl}\\[2mm]
	Nicolaus Copernicus University, Poland\\[6mm]
 {\em Dedicated to Professor Vygantas Paulauskas}}
	\date{}
	
	\maketitle
	
	\begin{abstract}
			For a stationary sequence that is regularly varying and associated we give conditions which guarantee that partial sums of this sequence, 
		under normalization related to the exponent of regular variation, converge in distribution to a stable, 
		non-Gaussian limit. The obtained limit theorem admits a natural extension to the functional convergence in Skorokhod's $M_1$ topology.
	\end{abstract}
	
	
	\noindent{\em MSClassification 2010:} 60F05, 60F17, 60E07, 60E15, 60G10.

\section{Introduction}\label{s:1}

In the pioneering work \cite{Pau76} Paulauskas considered some covariance-like quantities, which were defined for jointly stable random variables. He returned to this topic in recent papers \cite{DaPa14} and \cite{DaPa17} with a variety of examples based on linear processes with heavy tailed innovations (see also \cite{Pau13} and \cite{Pau19}). It is clear by now that these various quantities are very useful in limit theorems.

Perhaps most elegant example of limit theorems operating with covariance-like quantities was given in \cite{DaJa94} for associated and jointly stable stationary sequences (see also \cite{DaPa17} for generalization of these results to stationary random fields).

Recall (see e.g. \cite{SaTa94}) that random variables $X_1, X_2, \ldots $ are  jointly $\alpha$-stable, $0
< \alpha < 2$, if for each $n \in \GN$ there exists a finite Borel measure
$\Gamma_n$ on the unit sphere in $\GR^n$,
\[
\cS^{n-1} =\{ {\bf s}=(s_1,\ldots s_n)\in \GR^n\ :\
\sum_{i=1}^{n}s_{i}^{2}=1\}
\]
and a vector ${\bf
	b}_n \in \GR^n$ such that the characteristic function of ${\bf X}_n =
(X_1, X_2, \ldots, X_n)$ is of the form

\bequ{e4m}
E \exp i({\bf t},{\bf X}_n) =
\exp \left( i({\bf b}_n,{\bf t}) +  \int_{\cS^{n-1}}
\int_0^{\infty} g({\bf t},{\bf s},r) \frac{dr}{r^{\alpha+1}}
\Gamma_n(d\/{\bf s})   \right).
\eequ
Here 
\bequ{e4}
g({\bf t},{\bf s},r) = \left\{
\begin{array}{ll}
	e^{i ({\bf t},{\bf s}) r} - 1 ,\  &\mbox{\rm if } 0 < \alpha < 1, \\
	e^{i ({\bf t},{\bf s}) r} - 1 - i ({\bf t},{\bf s}) r I(r \leq 1),
	& \mbox{\rm if } \alpha = 1, \\
	e^{i ({\bf t},{\bf s}) r} - 1 - i ({\bf t},{\bf s})  r,\
	&\mbox{\rm if } 1 < \alpha  < 2.
\end{array}
\right.
\eequ
We will write
$ \cL({\bf X}_n) = \gamma_{\alpha}({\bf b}_n, \Gamma_n)$.
Clearly, if $X_1, X_2, \ldots $ are also strictly
stationary, then for some $b \in \GR^1$
\bequ{e5}
{\bf b}_n = ( \underbrace{b,\ldots, b}_{\mbox{\small $n$ times}} ),\quad
n\in \GN.
\eequ
For one-dimensional stable distributions we have
$\cS^0 = \{ -1, 1\}$, and we shall use  the notation
\[   \gamma_{\alpha}(b, \Gamma(\{1\}), \Gamma(\{-1\})) \equiv \gamma_{\alpha}(b,
\Gamma). \] 

Jointly stable random variables $X_1, X_2, \ldots$ are {\em strictly}
$\alpha$-stable, if either
\newline (a) ${\bf b}_n = 0,\ n \in \GN$, when $\alpha \neq 1$, or 
\newline (b) $\int_{\cS^{n-1}} {\bf s} \Gamma_n(d{\bf s}) = 0,\ n \in \GN$,
when $\alpha=1$. This holds if, for instance, $\Gamma_n$ is a symmetric measure
on $\cS^{n-1}$. 

The other key assumption in this paper is {\em association}.
Following \cite{EPW67} we call random variables $X_1, X_2, \ldots, X_n$  associated if
\begin{equation}
\Cov{f(X_1,X_2,\ldots,X_n)}{g(X_1,X_2,\ldots,X_n)} \geq 0,
\end{equation}
for each pair of functions $f, g : \GR^n \to \GR^1$, which are
non-decreasing in each coordinate and for which the above covariance
exists. An infinite collection  of random variables is associated, if its
every finite subset consists of associated random variables. 
For basic properties of associated random variables we refer to the original article \cite{EPW67} and also to the more recent source \cite{BuSh07}.

For jointly $\alpha$-stable random variables there exists an astonishingly simple description of association, due to Lee, Rachev and Samorodnitsky \cite{LRS90} (see also \cite{Res88}, \cite{Sam95} and \cite{HPS98} for related results). 
The measure $\Gamma_n$  
has to be concentrated on ``positive'' and ``negative'' parts of
$\cS^{n-1}$, i.e.
\bequ{e6}
\Gamma_n \Big( \cS^{n-1} \cap \Big\{ [0,+\infty)^n \cup (-\infty,0]^n
\Big\}^c \Big) = 0.
\eequ
This property allows proving  very nice limit theorems, which we restate here from \cite{DaJa94} as Theorems 1-3, for they will be used in the course of the proofs of our new, more general result.

In what follows
$S_n$ will always stand for $X_1 + X_2 + \ldots + X_n$. 

As in the case of independent summands,
we have separate results for the three cases where
$0 < \alpha < 1$, $\alpha=1$, and $1< \alpha <2$.

\begin{theorem}\label{th2.1}
	Let $X_1, X_2, \ldots $ be stationary, associated and jointly $\alpha$-stable,
	$0 < \alpha < 1$. Then 
	\bequ{e7}
	\frac{S_n}{n^{1/\alpha}} \indist \mu_{\infty},
	\eequ
	where $\mu_{\infty}$ is a strictly $\alpha$-stable distribution.
\end{theorem}

\begin{theorem}\label{th2.2}
	Let $X_1, X_2, \ldots$ be stationary, associated and jointly $1$-stable.
	Then there exist constants $A_n$ such that
	\bequ{e8}
	\frac{S_n}{n} - A_n \sim X_1.
	\eequ
	In particular, if $\Gamma_n$ is symmetric for each $n\in\GN$, then  
	\bequ{e9}
	\frac{S_n}{n} \sim X_1,\quad n\in\GN.
	\eequ 
\end{theorem}

\begin{theorem}\label{th2.3}
	Let $X_1, X_2, \ldots $ be stationary, associated and jointly $\alpha$-stable,
	$1 < \alpha < 2$, with two-dimensional distributions $
	\cL\Big((X_1,X_k)\Big) = \gamma_{\alpha}\Big((b,b), \Gamma_{\{1,k\}}\Big)$.
	
	If
	\bequ{e11}
	\sum_{k=2}^{\infty} \int_{\cS^1} s_1 s_2\, \Gamma_{\{1,k\}} (d{\bf s})  <
	+\infty, 
	\eequ
	then
	\bequ{e11a}
	\frac{S_n - E S_n}{n^{1/\alpha}} = \frac{S_n - n b}{n^{1/\alpha}} \indist \mu_{\infty},
	\eequ
	where $\mu_{\infty}$ is a non-degenerate strictly $\alpha$-stable distribution. 
\end{theorem}

Notice that Theorem \ref{th2.3} convincingly supports the point of view of \cite{Pau76} that the spectral covariance 
\[\int_{\cS^1} s_1 s_2\, \Gamma_{\{1,k\}} (d{\bf s})\] is a covariance-like quantity. To see this let us compare the shape of Theorem  \ref{th2.3} with the classic central limit theorem for associated stationary sequences due to Newman \cite{New80} (see also \cite{NeWr81} for functional convergence and \cite{BDD86} for a more general result): 
suppose  that $\{X_j\}$ is stationary and associated and that $E X_1 = 0$, $E X_1^2 < +\infty$. Then condition
\bequ{e2}
\sigma^2 = E X_1^2 + 2 \sum_{j=2}^{\infty} E X_1 X_j < +\infty,
\eequ
implies $
S_n/\sqrt{n} \indist \cN(0,\sigma^2)$, as $n\to\infty$.

Of course the assumption on joint $\alpha$-stability appearing in Theorems \ref{th2.1}, \ref{th2.2} and \ref{th2.3} is very restrictive. Theorem 2.8 in \cite{DaJa94} gets rid of this limitation and introduces another covariance-like quantity.

For a strictly stationary sequence
$\{X_j\}_{j\in\GN}$ define
\bequ{e3.7}
H_{(X_i,X_j)}(x_i,x_j) = P(X_i\leq x_i, X_j\leq x_j)
- P(X_i\leq x_i) P(X_j\leq x_j).
\eequ
Then fix $A>0$ and $\alpha \in (0,2)$ and define
\bequ{e3.11}
I^A_{\alpha}(X_i,X_j) = \sup_{a \geq A} a^{p-2} \int_{-a}^{a} \int_{-a}^{a}
H_{(X_i,X_j)}(x,y)\, dx\, dy.
\eequ
It is immediate that if $\{X_j\}$ is associated, then both  $H_{(X_i,X_j)}(x_i,x_j) \geq 0$ and $I^A_{\alpha}(X_i,X_j) \geq 0$, with the latter taking possibly the value $+\infty$. The quantity $I^A_{\alpha}(X_i,X_j)$ satisfies the Cauchy-Schwarz inequality:
\[
I^A_{\alpha}(X_i,X_j) \leq \sqrt{I_{\alpha}^A(X_i,X_i)}\sqrt{I_{\alpha}^A(X_j,X_j)}.
\]
Moreover, as \cite[Theorem 2.8]{DaJa94} states, condition
\bequ{e3.12}
\sum_{k=2}^{\infty} I^A_{\alpha}(X_1,X_k) < +\infty,
\eequ
{\em plus} some natural distributional conditions imply convergence of partial sums to stable laws.

The serious drawback of  coefficient $I^A_{\alpha}(X_i,X_j)$ is that it is infinite for some  marginal distributions $\cL(X_i)$  belonging to the domain of attraction of a stable law.

In the present paper we give assumptions which are the most general when considering the framework based on domains of attraction and Newman's inequality.

\section{Statement of results}\label{s:2}

Let $\{X_j\}$ be a stationary sequence. We will say that it is 
{\em regularly varying}, if it is {\em jointly regularly varying} with some index $\alpha$, i.e. for each $i \leq j$ the joint distribution of $(X_i, X_{i+1}, \ldots, X_j)$ is regularly varying with necessarily the same index $\alpha$. We will use a reformulation of regular variation, which is close in spirit to \cite[Theorem 2.1]{BaSe09} and can be proved in a similar way (see also \cite{Res86}).

In what follows we shall assume that
\begin{equation}\label{eq:X0}
\bP(|X_1| > x) = x^{-\alpha} \ell(x), 
\end{equation}
where $\alpha \in (0,2)$ and $\ell(x)$ is a slowly varying function. Given \nref{eq:X0} we define the normalizing constants $B_n$ by the relation 
\begin{equation}\label{eq:been}
n \bP( |X_1| > B_n) \to 1. 
\end{equation}
It is well-known that $\{B_n\}$ is $1/\alpha$-regularly varying.

For multidimensional distributions we assume that for each $N\in \GN$
\begin{equation}
n \bP\Big( \big(\frac{X_1}{B_n},\frac{X_2}{B_n},\ldots, \frac{X_N}{B_n}\big) \in \big(\cdot\big)\Big) \conv \nu_N\big(\cdot\big),\ \text{ vaguely on $\overline{\GR}^m \setminus \{\mathbf{0}\}$},
\end{equation}
where $\nu_N$ is (necessarily) a L\'evy measure on $\GR^N$.

In addition we always assume that 
\begin{align} \label{eq:always1}
	\text{ if $\alpha = 1$,}&\ \text{ then vectors $(X_1, X_2, \ldots, X_N)$ have symmetric distributions, $N\in \GN$}, \\  
	\text{ if $\alpha \in (1,2)$,} &\ \text{ then $\bE X_1  = 0$}.\label{eq:always2}
\end{align}

Let us observe that our assumptions \nref{eq:X0} - \nref{eq:always2} imply that for each $N \in \GN$
\bequ{eq:3.1A}
\frac{{\bf Z}_{N,1} + {\bf Z}_{N,2} + \ldots + {\bf Z}_{N,n}}{B_n}
\indist (Y_1^N, Y_2^N, \ldots, Y_N^N), %\asn,
\eequ
where ${\bf Z}_{N,1}, {\bf Z}_{N,2}, \ldots, {\bf Z}_{N,n}, \ldots$ are
independent copies of $(X_1, X_2, \ldots, X_N)$ and $(Y_1^N, \ldots, Y_N^N)$ is a strictly  $\alpha$-stable random vector with the distribution determined by $\nu_N$.

For fixed $a > 0$ we define a function
$f_a : \GR^1 \to \GR^1$ by
\bequ{e3.9}
f_{a}(x)=\left\{\begin{array}{rl}
	a &\mbox{\rm if } x>a \\
	x &\mbox{\rm if } |x|\leq a \\
	-a &\mbox{\rm if } x < -a
\end{array} \right. .
\eequ
Note that $f_{a}(x/b)=b^{-1}f_{ab}(x)$, that
$f_{a}(x)$ is a non-decreasing function in $x$, and that
$\left\{ f_{a}(X_{j})\ :\ j\geq 1\right\}$ is again an associated
sequence of random variables. Moreover, $f_a(x)$ is absolutely continuous with $f_a'(x) = I_{(-a,a)}(x)$ a.e. and so, by \cite[Lemma 3.1]{Yu93} 
\bequ{e3.10}
\Cov{f_a(X_i)}{f_a(X_j)} = \int_{-a}^{a} \int_{-a}^{a}
H_{(X_i,X_j)}(x_i,x_j)\, dx_i\, dx_j.
\eequ
It follows from \nref{eq:3.1A} that for each $N \geq 2$ and every $a > 0$
\begin{equation}\label{eq:single}
\begin{split}n \cdot B_n^{-2}  \text{Cov}
\Big( f_{a\cdot B_n}(X_1), f_{a\cdot B_n}(X_N)\Big) & = n\cdot \text{Cov} \Big( f_a\big(\frac{X_1}{B_n}\big), f_a\big(\frac{X_N}{B_n}\Big) \\
& \mbox{}\qquad\qquad\conv_{n\to\infty}  \int_{\GR^2}
f_a(x_1) f_a(x_2)\,
\nu_{\{1,N\}}(dx_1, dx_2),
\end{split}
\end{equation}
where $\nu_{\{1,N\}}(dx_1, dx_2)$ is the L\'evy measure of 
$\cL(Y_1^N,Y_N^N)$ (see e.g. \cite[Theorem 2.35, p.362]{JaSh87}).
This implies that for each $N \geq 2$ the non-decreasing function
\[ a \mapsto g_N(a) = \text{Cov} \big(f_a(X_1), f_a(X_N)\big) \]
is regularly varying and the exponent of regular variation is $2-\alpha$ (see e.g. \cite[Lemma 3, p. 277]{Fel71}).

In fact our main assumption requires substantially more:
\begin{equation}  \label{eq:more}
a \mapsto \sum_{j=2}^{\infty} g_j(a)\ \  \text{ is a regularly varying function}.
\end{equation}
We relate this abstract property with \nref{eq:single} by assuming 
\begin{equation}\label{eq:part}
n \sum_{j=2}^{\infty} 
{\rm Cov}\Big(f_{1}\big(\frac{X_1}{B_n}\big),f_{1}\big(\frac{X_j}{B_n}\big)\Big) \conv_{n\to\infty}  \sum_{j=2}^{\infty} \int_{\GR^2}
f_1(x_1) f_1(x_2)\,
\nu_{\{1,j\}}(dx_1, dx_2) < +\infty.
\end{equation}
Both \nref{eq:more} and \nref{eq:part} imply that $\sum_{j=2}^{\infty} g_j(a)$ is $(2-\alpha)$-regularly varying.

\begin{remark} Assumptions \nref{eq:more} and \nref{eq:part} taken together are equivalent to 
	\begin{equation}\label{eq:partduo}
	n \sum_{j=2}^{\infty} 
	{\rm Cov}\Big(f_{a}\big(\frac{X_1}{B_n}\big),f_{a}\big(\frac{X_j}{B_n}\big)\Big) \conv_{n\to\infty}  \sum_{j=2}^{\infty} \int_{\GR^2}
	f_a(x_1) f_a(x_2)\,
	\nu_{\{1,j\}}(dx_1, dx_2) < +\infty, \ \ a > 0.
	\end{equation}
	As the limit is continuous and monotone in $a$, the above pointwise convergence is, in fact, uniform on bounded intervals.
\end{remark}

\begin{theorem}\label{th:main}
	Let $\{X_j\}$ be a stationary sequence that is associated and satisfies conditions \nref{eq:X0} - \nref{eq:always2}. 
	
	Suppose that \nref{eq:more} and \nref{eq:part} hold with $B_n$ defined by \nref{eq:been}. 
	
	Then there exists a strictly $\alpha$-stable distribution $\mu_{\infty}$ such that
	\[ \frac{X_1 + X_2 + \ldots + X_{n}}{B_n} \indist \mu_{\infty}.\]  
\end{theorem}

\begin{remark} Theorem 2.13 in \cite{DaJa94} operates with apparently weaker assumption on domain of attraction for sums of $S_N$, $N\in\GN$, only. It is not clear whether assumptions of this type are weaker or equivalent to our conditions \nref{eq:X0} - \nref{eq:always2}.  See \cite{BDM02} for discussion of problems of similar flavor.
\end{remark}

\begin{remark}
	As Remark 2.4 in \cite{DaJa94} shows, if $\alpha \in (0,1)$ then there are jointly stable associated sequences such that the limit is {\em degenerate}. We think that in the presence of the strong assumption of summability \nref{eq:partduo} it is possible to show the non-degeneracy of the limit. We are, however, not able to prove this statement.  
\end{remark}

Theorem \ref{th:main} admits a natural functional extension.
Let us define a sequence of stochastic processes with trajectories in the Skorokhod space $\GD\big([0,1]:\GR^1\big)$.
\begin{equation}\label{eq:process}
S_n(t) = \frac{S_{\lfloor n\cdot t\rfloor}}{B_n},\ \ t\in [0,1].
\end{equation}
It is known (\cite{AvTa92}) that in the general setting of associated sequences it is impossible to obtain the convergence of $\{S_n(t)\}$ in Skorokhod's $J_1$ topology. On the other hand, applying the powerful Theorem 1 of \cite{LoRi11}, we shall obtain the convergence in Skorokhod's $M_1$ topology. For the definitions and basic properties of Skorokhod's topologies we refer either to the seminal paper \cite{Skor56} or to the extensive source \cite{Whit02}.

\begin{theorem}\label{th:funct} In assumptions of Theorem \ref{th:main}, 
	the sequence $\{S_n(t)\}$ converges in law on the Skorokhod space $\GD\big([0,1]:\GR^1\big)$ equipped with Skorokhod's $M_1$ topology. The limit is the stable L\'evy process $\{Y(t)\}$ given by $Y(1) \sim \mu_{\infty}$.
\end{theorem}

\section{Proofs}
As noted in Introduction, we  follow the line of the proof of Theorem 2.8 in \cite{DaJa94}. But the details are different  in many places, for our result is more general. Therefore we give 
here a complete proof.

\subsection{Proof of Theorem \ref{th:main}}

Let us recall \nref{eq:3.1A}, i.e. for each $N$

\[
\frac{{\bf Z}_{N,1} + {\bf Z}_{N,2} + \ldots + {\bf Z}_{N,n}}{B_n}
\indist (Y_1^N, Y_2^N, \ldots, Y_N^N), %\asn,
\]
where ${\bf Z}_{N,1}, {\bf Z}_{N,2}, \ldots, {\bf Z}_{N,n}, \ldots$ are
independent copies of $(X_1, X_2, \ldots, X_N)$. It follows that there exists a stationary process $\{Y_j\}$ such that
\[ (Y_1,Y_2,\ldots, Y_N) \sim (Y_1^N, Y_2^N, \ldots, Y_N^N), \ N\in \GN.
\]
The process $\{Y_j\}$ is {\em jointly $\alpha$-strictly stable} and {\em associated}. We shall call $\{Y_j\}$ the stable tangent to $\{X_j\}$, for the asymptotic properties of partial sums of the original and the tangent processes are the same. Notice that for $\alpha \in [1,2)$ the existence of the tangent process requires more than just the regular variation of $\{X_j\}$, therefore we do not introduce here the {\em tail} process as defined in \cite{BaSe09}.

Let $\mu_N$ be the distribution of $Y_1 + Y_2 + \ldots + Y_N$.
By the strict $\alpha$-stability 
\begin{equation}\label{eq:3.2}
\frac{Y_1 + Y_2 + \ldots + Y_N}{N^{1/\alpha}} \sim \mu_N^{*(1/N)}.
\end{equation} 
Here $\mu^{*\beta}$ is the convolution $\beta$-power of the infinitely divisible distribution $\mu$.

By the association and Theorems \ref{th2.1}, \ref{th2.2} and  \ref{th2.3}
there exists a strictly $\alpha$-stable distribution $\mu_{\infty}$ such that as $N\to\infty$
\begin{equation}\label{eq:3.3}
\mu_N^{*(1/N)} \wconv \mu_{\infty}.
%\frac{Y_1 + Y_2 + \ldots + Y_N}{N^{1/\alpha}} \indist 
\end{equation}
Notice that for $\alpha \in (1,2)$ relation \nref{e11} in Theorem \ref{th2.3} is satisfied by \cite[Remark 2.6]{DaJa94} and our assumption \nref{eq:part}.

It follows from \nref{eq:3.1A} and \nref{eq:3.2} that for each $N\in \GN$ and as $n\to \infty$
\[ \big(\bE e^{i\lambda (X_1 + X_2 + \ldots + X_N)/ N^{1/\alpha} B_n}\big)^n \conv \big(\hat{\mu}_N(\lambda)\big)^{1/N}, \ \lambda \in \GR^1.\]
By the regular variation of $\{B_n\}$, $N^{1/\alpha}B_n \sim B_{N\cdot n}$, as $n\to \infty$, so the above relation can be rewritten as 
\[ \big(\bE e^{i\lambda S_N/ B_n)}\big)^{\lfloor n/N\rfloor} \conv \big(\hat{\mu}_N(\lambda)\big)^{1/N}, \ N\in \GN, \ \lambda \in \GR^1.\]
This and \nref{eq:3.3} imply 
$
\lim_{N \to \infty} \lim_{n\to \infty} \big| \big(\bE e^{i\lambda S_N/ B_n)}\big)^{\lfloor n/N\rfloor} -  \hat{\mu}_{\infty}(\lambda)\big| = 0, \ \lambda \in \GR^1.
$
Therefore it is enough to prove that $
\lim_{N\to\infty}
\limsup_{n\to\infty} \big| \bE e^{i \lambda (S_n/B_n)} -
\big( E e^{i \lambda 
	(S_N/B_n)}\big)^{\lfloor n/N\rfloor} \big| = 0, \ \lambda \in \GR^1,
$ or, after a simple modification,
\begin{equation}\label{eq:3.4}
\lim_{N\to\infty}
\limsup_{m\to\infty} \Big| \bE e^{i \lambda S_{m \cdot N}/B_{m\cdot N}} -
\Big( E e^{i \lambda 
	S_N/B_{m \cdot N}}\Big)^m \Big| = 0,\ \ \lambda\in \GR^1.
\end{equation}

Recall that function $f_a$ is defined by \nref{e3.9}.
Consider the following decomposition.
\begin{equation*}
	B_n^{-1} \sum_{j=1}^{k}X_{j}
	=   \sum_{j=1}^{k}f_{a}(B_n^{-1}X_{j} )
	+ \sum_{j=1}^{k}\Big( B_n^{-1} X_{j}
	- f_{a}(B_n^{-1}X_{j})\Big) 
	=:  T^{(a)}_{n,k}+V_{n,k}^{(a)}.
\end{equation*}
Choose arbitrary $\eta > 0$. We have for $a > \eta^{-1/\alpha}$ 
\[
\begin{array}{rcl}
\limsup_{n\to\infty} \bP(V_{n,k}^{(a)}\neq 0) & \leq & \limsup_{n\to\infty} \bP\Big( \exists\ 1\leq j\leq n\ :\
|X_{j}|>a B_n\Big) \\
& \leq & \limsup_{n\to\infty} n\bP\Big( |X_{1}|>a B_n\Big) =   a^{-\alpha} < \eta. 
\end{array}
\]
Consequently \[
\limsup_{n\to\infty} \Big| E e^{i \lambda S_n/B_n}
- E e^{i \lambda T_{n,n}^{(a)}} \Big| < 2\eta.\] 
A similar reasoning also shows that
\[
\limsup_{m\to\infty}\Big| \Big(E e^{i \lambda 
	S_N/B_{m\cdot N}}\Big)^m - \Big(E e^{i
	\lambda  T_{m\cdot N, N}^{(a)} }\Big)^m \Big| < 2\eta.
\]
It follows that \nref{eq:3.4} will hold provided for each $a > 0$
\begin{equation}\label{eq:3.4A}
\lim_{N\to\infty}
\limsup_{m\to\infty} \Big| E e^{i \lambda \sum_{j=1}^{m\cdot N} U_{m\cdot N,j}^{(a)}} - \Big(E e^{i
	\lambda \sum_{j=1}^{N} U_{m\cdot N,j}^{(a)}  }\Big)^m \Big| = 0,
\end{equation}
where 
\[ U_{n,j}^{(a)} = f_a(B_n^{-1} X_j).\]
Now we are ready to apply Newman's inequality \cite{New80} (see also \cite[Theorem 1]{NeWr81}). Take $\lambda \in \GR^1$ and return for a while to $n=m\cdot N$. Then
\begin{align*}
	\Big| 
	E \exp &\Big\{ i\lambda\sum_{j=1}^{n} U_{n,j}^{(a)}\Big\} -  
	\Big( E \exp \Big\{ i\lambda\sum_{j=1}^{N}
	U_{n,j}^{(a)}\Big\}\Big)^{m}
	\Big| \\
	& \leq \frac{\lambda^2}{2} \sum_{1 \leq k \neq l \leq m} \text{Cov}\Big(\sum_{i=(k-1)\cdot N + 1}^{k\cdot N} U_{n,i}^{(a)}\Big)\Big(\sum_{j=(l-1)\cdot N + 1}^{l\cdot N} U_{n,j}^{(a)}\Big)\qquad \text{Newman's inequality}
	\\
	&= \frac{\lambda^2}{2} \Big( \text{Var}\Big( \sum_{j=1}^{n} U_{n,j}^{(a)}\Big) - m  \text{Var}\Big( \sum_{j=1}^{N} U_{n,j}^{(a)}\Big)\Big) \\
	&=
	\frac{\lambda^2 n}{2} \Big( \frac{1}{n} \text{Var}\Big( \sum_{j=1}^{n} U_{n,j}^{(a)}\Big) - \frac{1}{N} \text{Var}\Big( \sum_{j=1}^{N} U_{n,j}^{(a)}\Big)\Big)\\
	&=  \lambda^{2}\Bigg\{
	\sum_{j=2}^{N}\Big(\frac{1}{N}-\frac{1}{n}\Big)(j-1)
	\Big( n \Cov{U_{n,1}^{(a)}}{U_{n,j}^{(a)}}\Big) \\
	&\qquad\qquad\qquad  + \sum_{j=N+1}^{n}\Big(1-\frac{j-1}{n}\Big) 
	\Big( n \Cov{U_{n,1}^{(a)}}{U_{n,j}^{(a)}}\Big)
	\Bigg\} \\
	&\leq  \lambda^{2}\left\{ \frac{1}{N} \sum_{i=1}^{N}\sum_{j=i+1}^N \Big(n
	\Cov{U_{n,1}^{(a)}}{U_{n,j}^{(a)}} \Big)  + \frac{1}{N} \sum_{i=1}^{N}\sum_{j=N+1}^n \Big(n
	\Cov{U_{n,1}^{(a)}}{U_{n,j}^{(a)}} \Big) \right\} \\
	&= \lambda^{2} \frac{1}{N} \sum_{i=1}^N \sum_{j=i+1}^n \Big( n
	\Cov{U_{n,1}^{(a)}}{U_{n,j}^{(a)}} \Big) \\
	&\leq \lambda^{2} \frac{1}{N} \sum_{i=1}^N n B_n^{-2} \sum_{j=i+1}^{\infty} 
	{\rm Cov}\Big(f_{a\cdot B_n}(X_1),f_{a\cdot B_n}(X_j)\Big) \\
	&\conv_{m\to \infty} \lambda^{2} a^{2-\alpha} \frac{1}{N} \sum_{i=1}^N  \sum_{j=i+1}^{\infty} \int_{\GR^2}
	f_1(x_1) f_1(x_2)\,
	\nu_{\{1,j\}}(dx_1, dx_2) \conv_{N\to \infty} 0,
\end{align*}
by \nref{eq:single}, \nref{eq:more} and \nref{eq:part}.
\subsection{Proof of Theorem \ref{th:funct}}

In view of \cite[Theorem 1]{LoRi11} it is enough to establish the finite dimensional convergence. Since the increments of $\{S_n(t)\}$ are (asymptotically) stationary we need only asymptotic independence of the increments. For the sake of brevity we shall restrict our attention to two adjoining increments  $S_n(t_1) - S_n(t_0)$ and $S_n(t_2) - S_n(t_1)$, $ 0 \leq t_0 < t_1 < t_2 \leq 1$. Let $\lambda, \theta \in \GR^1$. 

In order to prove that 
\begin{align*} \lim_{n\to\infty} \Big|\bE \exp \Big( i\lambda \big( S_n(t_1) &- S_n(t_0)\big) + i\theta \big( S_n(t_2) - S_n(t_1)\big) \Big) \\
	&-  \bE \exp \Big( i\lambda \big( S_n(t_1) - S_n(t_0)\big)\Big) \cdot \bE \exp \Big( i\theta \big( S_n(t_2) - S_n(t_1)\big) \Big) \Big| =0,
\end{align*}
we may, as before replace the increments with sums of $U_{n,j}^{(a)}$, for $a > 0$ large enough.
Then we have
\begin{align*}
	\Big|\bE \exp \Big\{ &i\lambda \Big(\sum_{ \lfloor n t_0\rfloor< j \leq \lfloor n t_1\rfloor} U_{n,j}^{(a)}\Big) + i\theta \Big(\sum_{ \lfloor n t_1\rfloor< j \leq \lfloor n t_2\rfloor} U_{n,j}^{(a)}\Big) \Big\} \\
	&-  \bE \exp \Big\{ i\lambda \Big(\sum_{ \lfloor n t_0\rfloor< j \leq \lfloor n t_1\rfloor} U_{n,j}^{(a)}\Big) \Big\} \cdot \bE \exp \Big\{ i\theta \Big(\sum_{ \lfloor n t_1\rfloor< j \leq \lfloor n t_2\rfloor} U_{n,j}^{(a)}\Big) \Big\} \Big| \\
	& \leq |\lambda||\theta|\text{Cov} \Big(\sum_{ \lfloor n t_0\rfloor< j \leq \lfloor n t_1\rfloor} U_{n,j}^{(a)}, \sum_{ \lfloor n t_1\rfloor< k \leq \lfloor n t_2\rfloor} U_{n,j}^{(a)}\Big)\qquad\qquad \text{Newman's inequality} \\
	&\leq |\lambda||\theta| \sum_{r=1}^{n-1} r\cdot \text{Cov} \big(U_{n,1}^{(a)},U_{n,1 + r}^{(a)}\big) \\
	&\leq  |\lambda||\theta|\Bigg\{
	\frac{M}{n} \sum_{r=1}^M  n B_n^{-2}
	{\rm Cov}\Big(f_{a\cdot B_n}(X_1),f_{a\cdot B_n}(X_{1+r})\Big)\\
	&\mbox{}\qquad\qquad\qquad + \sum_{r=M+1}^{\infty} n B_n^{-2}
	{\rm Cov}\Big(f_{a\cdot B_n}(X_1),f_{a\cdot B_n}(X_{1+r})\Big) \Bigg\} \\
	&\conv_{n\to\infty}  |\lambda||\theta| a^{2-\alpha}\sum_{r=M+1}^{\infty} \int_{\GR^2}
	f_1(x_1) f_1(x_2)\,
	\nu_{\{1,1+r\}}(dx_1, dx_2) \conv_{M\to \infty} 0,
\end{align*}
 again by \nref{eq:single}, \nref{eq:more} and \nref{eq:part}.

\end{document}